\documentclass[12pt,reqno]{article}

\usepackage[usenames]{color}
\usepackage{amssymb}
\usepackage{graphicx}
\usepackage{amscd}

\usepackage[colorlinks=true,
linkcolor=webgreen,
filecolor=webbrown,
citecolor=webgreen]{hyperref}

\definecolor{webgreen}{rgb}{0,.5,0}
\definecolor{webbrown}{rgb}{.6,0,0}

\usepackage{color}
\usepackage{fullpage}
\usepackage{float}

\usepackage{graphics,amsmath,amssymb}
\usepackage{amsthm}
\usepackage{amsfonts}
\usepackage{latexsym}
\usepackage{epsf}

\newcommand{\seqnum}[1]{\href{http://www.research.att.com/cgi-bin/access.cgi/as/~njas/sequences/eisA.cgi?Anum=#1}{\underline{#1}}}

\begin{document}

%\begin{center}
%\epsfxsize=4in
%\leavevmode\epsffile{logo129.eps}
%\end{center}

\begin{center}
\vskip 1cm{\LARGE\bf
A Combinatorial Interpretation for  \\
\vskip .1in
 Certain Relatives of the Conolly Sequence
}
\vskip 1cm \large
B. Balamohan, Zhiqiang Li, and Stephen Tanny\footnote{Corresponding author.}\\
Department of Mathematics\\
University of Toronto\\
Toronto, Ontario M5S 2E4\\
Canada\\
\href{mailto:bbalamoh@math.utoronto.ca}{\tt bbalamoh@math.utoronto.ca}\\
\href{mailto:zhiqiang.li@utoronto.ca}{\tt zhiqiang.li@utoronto.ca}\\
\href{mailto:tanny@math.utoronto.ca}{\tt tanny@math.utoronto.ca}\\
\end{center}

\vskip .2 in

\begin{abstract}
For any integer $s\geq0$, we derive a combinatorial interpretation
for the family of sequences generated by the recursion
(parameterized by $s$)
$h_s(n)=h_s(n-s-h_s(n-1))+h_s(n-2-s-h_s(n-3)), n > s+3,$ with the
initial conditions $h_s(1) = h_s(2) = \cdots = h_s(s+2) = 1$ and
$h_s(s+3) = 2$. We show how these sequences count the number of
leaves of a certain infinite tree structure. Using this
interpretation we prove that $h_{s}$ sequences are ``slowly
growing'', that is, $h_{s}$ sequences are monotone nondecreasing,
with successive terms increasing by 0 or 1, so each sequence hits
every positive integer. Further, for fixed $s$ the sequence $h_s(n)$
hits every positive integer twice except for powers of 2, all of
which are hit $s+2$ times. Our combinatorial interpretation provides
a simple approach for deriving the ordinary generating functions for
these sequences.

\end{abstract}

\theoremstyle{plain}
\newtheorem{theorem}{Theorem}[section]
\newtheorem{corollary}[theorem]{Corollary}
\newtheorem{conjecture}[theorem]{Conjecture}
\newtheorem{question}[theorem]{Question}
\newtheorem{lemma}[theorem]{Lemma}
\newtheorem{proposition}[theorem]{Proposition}
\theoremstyle{definition}
\newtheorem{definition}[theorem]{Definition}
\theoremstyle{remark}
\newtheorem*{notation}{Notation}
\newtheorem*{remark}{Remark}

\numberwithin{equation}{section}
\numberwithin{theorem}{section}
\numberwithin{table}{section}
\numberwithin{figure}{section}

\newcommand{\ve}{\varepsilon}
\newcommand{\C}{\mathbb C}
\newcommand{\N}{\mathbb N}
\newcommand{\R}{\mathbb R}
\newcommand{\Z}{\mathbb Z}
\newcommand{\Q}{\mathbb Q}
\renewcommand{\H}{\mathcal H}

\def\({\left(}
\def\){\right)}
\def\[{\left[}
\def\]{\right]}

\def \mm #1#2#3#4{\begin{pmatrix} #1 & #2 \cr #3 & #4 \end{pmatrix}}

\newcommand{\ontop}[2]{\genfrac{}{}{0pt}{}{#1}{#2}}

\section{Introduction}
Conolly \cite{ConVa} introduced the sequence defined by the
meta-Fibonacci (self-referencing) recursion
\begin{equation}\label{equConF}
F(n) = F(n - F(n-1)) + F(n - 1 - F(n-2)),\qquad       n>2
\end{equation}
with initial conditions $(F(1) = 0, F(2) = 1)$ or $(F(1) = 1,
F(2) = 1)$. He notes that for $n>2$ the recursion yields the same
sequence whether $F(1) = 0$ or $F(1) = 1$. The resulting sequence
behaves in a very simple fashion: it is monotonic non-decreasing,
with successive terms differing by either $0$ or $1$, so the
sequence hits every positive integer. Following Ruskey \cite{FR} we
term such a sequence ``slowly growing''.

Variants and generalizations of the recursion (\ref{equConF}) and
its initial conditions have been studied (see, for example,
\cite{CCT}, \cite{CElz}, \cite{HiTan}, \cite{Tanny}). The most
comprehensive of these is the following $k$-term generalization,
where $s$ and $k \geq 2$ are nonnegative integers:
\begin{equation}\label{equTsk}
T_{s,k}(n) = \sum_{i=0}^{k-1} T_{s,k}(n-i-s-T_{s,k}(n-i-1)),\qquad
n>s+k.
\end{equation}

The special case of (\ref{equTsk}) for $s =1$ and initial conditions
$T_{1,k}(1) = T_{1,k}(2) = 1$, and $T_{1,k}(j) = j-1$ for $3\leq j
\leq k+1$ is first discussed by Higham and Tanny \cite{HiTan}. It is shown there that
the resulting sequence is slowly growing.

Recently Jackson and Ruskey \cite{BJaFRus} (and see \cite{CDegFRus}
for a generalization) derived a beautiful combinatorial
interpretation for the sequences generated by (\ref{equTsk}) with
initial conditions $T_{s,k}(j) = 1$ for $1 \leq j \leq s+1$, and
$T_{s,k}(j) = j-s$ for $s+2 \leq j \leq s+k$. These initial
conditions are natural extensions for arbitrary $s$ of those
appearing in \cite{HiTan}. Using their interpretation they show that
all of these sequences are slowly growing.

Even further, Jackson and Ruskey's interpretation provides a highly
intuitive understanding for the similarity of these sequences for
fixed $k$ and varying values of the ``shift'' parameter $s$, a fact
that had been noted by others but never proved. For fixed $k$ the
only differences between the sequences $T_{s,k}(n)$ and $T_{0,k}(n)$
occur in the frequencies with which each sequence hits the powers of
$k$. More precisely, for fixed $k$ and any positive integer $r$, the
value $k^r$ occurs $s$ more times in the sequence $T_{s,k}(n)$ than
in the sequence $T_{0,k}(n)$.

The present work is inspired by their approach. We begin with a
meta-Fibonacci recursion closely related to (\ref{equTsk}), namely,
\begin{equation}\label{equh_sk}
h_{s,k}(n) = \sum_{j=0}^{k-1}
h_{s,k}(n-2j-s-h_{s,k}(n-(2j+1))),\qquad n>s+2k+1.
\end{equation}

The right hand side of (\ref{equh_sk}) includes precisely those
terms in the sum on the right hand side of (\ref{equTsk}) that
derive from the even values of the index. For $k>2$ there does not
appear to be any set of initial conditions for which (\ref{equh_sk})
generates an interesting infinite sequence. For $k = 2$, however,
(\ref{equh_sk}) reduces to the two term recurrence
\begin{equation}     \label{eqnHs}
h_s(n)=h_s(n-s-h_s(n-1))+h_s(n-2-s-h_s(n-3)),\qquad  n > s+3.
\end{equation}
For any integer $s \geq 0$ we derive a combinatorial interpretation
for the sequence generated by the recursion (\ref{eqnHs}) with the
initial values $h_s(1) = h_s(2) = \cdots = h_s(s+2) = 1$ and
$h_s(s+3)=2$. These initial conditions are the natural analogue of
the ones described above by Deugau and Ruskey \cite{CDegFRus}. It turns out that this
sequence is slowly growing. Further, the role played by the
parameter $s$ in this family of sequences appears to be the same as
in the generalized Conolly recursion analyzed by Deugau and Ruskey \cite{CDegFRus}.
For given $s$ each sequence $h_s(n)$ hits every positive integer
twice except for the powers of 2, all of which are hit $(s+2)$ times
in the sequence. See Table \ref{tab1} for the first 50 terms of each
sequence for $0\leq s \leq 6$.

\begin{table}[!ht]
\fontsize{8}{8}\selectfont \caption{First 50 entries of $h_s$ for
$0\leq s \leq 6$.}\label{tab1} \center{
\begin{tabular}{ |c|c c c c c c c|}
\hline $n\backslash s$
     &   0   &   1   &   2   &   3   &   4   &   5   &   6   \\
\hline
1    &   1   &   1   &   1   &   1   &   1   &   1   &   1   \\

2    &   1   &   1   &   1   &   1   &   1   &   1   &   1   \\

3    &   2   &   1   &   1   &   1   &   1   &   1   &   1   \\

4    &   2   &   2   &   1   &   1   &   1   &   1   &   1   \\

5    &   3   &   2   &   2   &   1   &   1   &   1   &   1   \\

6    &   3   &   2   &   2   &   2   &   1   &   1   &   1   \\

7    &   4   &   3   &   2   &   2   &   2   &   1   &   1   \\

8    &   4   &   3   &   2   &   2   &   2   &   2   &   1   \\

9   &   5   &   4   &   3   &   2   &   2   &   2   &   2   \\

10   &   5   &   4   &   3   &   2   &   2   &   2   &   2   \\

11   &   6   &   4   &   4   &   3   &   2   &   2   &   2   \\

12   &   6   &   5   &   4   &   3   &   2   &   2   &   2   \\

13   &   7   &   5   &   4   &   4   &   3   &   2   &   2   \\

14   &   7   &   6   &   4   &   4   &   3   &   2   &   2   \\

15   &   8   &   6   &   5   &   4   &   4   &   3   &   2   \\

16   &   8   &   7   &   5   &   4   &   4   &   3   &   2   \\

17   &   9   &   7   &   6   &   4   &   4   &   4   &   3   \\

18   &   9   &   8   &   6   &   5   &   4   &   4   &   3   \\

19   &   10  &   8   &   7   &   5   &   4   &   4   &   4   \\

20   &   10  &   8   &   7   &   6   &   4   &   4   &   4   \\

21   &   11  &   9   &   8   &   6   &   5   &   4   &   4   \\

22   &   11  &   9   &   8   &   7   &   5   &   4   &   4   \\

23   &   12  &   10  &   8   &   7   &   6   &   4   &   4   \\

24   &   12  &   10  &   8   &   8   &   6   &   5   &   4   \\

25   &   13  &   11  &   9   &   8   &   7   &   5   &   4   \\

26   &   13  &   11  &   9   &   8   &   7   &   6   &   4   \\

27   &   14  &   12  &   10  &   8   &   8   &   6   &   5   \\

28   &   14  &   12  &   10  &   8   &   8   &   7   &   5   \\

29   &   15  &   13  &   11  &   9   &   8   &   7   &   6   \\

30   &   15  &   13  &   11  &   9   &   8   &   8   &   6   \\

31   &   16  &   14  &   12  &   10  &   8   &   8   &   7   \\

32   &   16  &   14  &   12  &   10  &   8   &   8   &   7   \\

33   &   17  &   15  &   13  &   11  &   9   &   8   &   8   \\

34   &   17  &   15  &   13  &   11  &   9   &   8   &   8   \\

35   &   18  &   16  &   14  &   12  &   10  &   8   &   8   \\

36   &   18  &   16  &   14  &   12  &   10  &   8   &   8   \\

37   &   19  &   16  &   15  &   13  &   11  &   9   &   8   \\

38   &   19  &   17  &   15  &   13  &   11  &   9  &   8   \\

39   &   20  &   17 &   16  &   14  &   12  &   10  &   8   \\

40   &   20  &   18  &   16  &   14  &   12  &   10  &   8   \\

41   &   21  &   18 &   16  &   15  &   13  &   11  &   9   \\

42   &   21  &   19  &   16  &   15  &   13  &   11  &   9   \\

43   &   22  &   19  &   17  &   16  &   14  &   12  &   10  \\

44   &   22  &   20  &   17  &   16  &   14  &   12  &   10  \\

45   &   23  &   20  &   18  &   16  &   15  &   13  &   11  \\

46   &   23  &   21  &  18   &   16  &   15  &   13  &   11  \\

47   &   24  &   21  &   19  &   16  &   16  &   14 &   12  \\

48   &   24  &   22  &   19  &   17  &   16  &   14  &   12  \\

49   &   25  &   22  &   20  &   17   &   16  &   15  &   13  \\

50   &   25  &   23  &   20  &   18   &   17  &   15  &   13   \\

\hline
\end{tabular}
}
\end{table}

In Section \ref{secTree} we show how the sequence
$\{h_s(n)\}_{n=1}^{\infty}$ can be interpreted in terms of the
leaves in an infinite tree structure, from which the properties
stated above follow readily. This interpretation also allows us to
derive the generating function for this sequence in Section
\ref{secGenFn}. We provide some final thoughts in Section
\ref{secRemark}.

\section{Tree Structure}  \label{secTree}
We begin this section by constructing a tree structure
$\mathcal{T}_{s}$ that consists of an infinite number of rooted trees
joined together at their respective roots plus a single initial
isolated node. Subsequently we show how to label the vertices in
$\mathcal{T}_{s}$ so that for any $s\geq 0$ the number of leaves up
to the vertex with label $n$ in this infinite tree structure is equal to $h_s(n)$.

\begin{figure}[htbp]
\begin{center}
 \epsfxsize=12cm\epsffile{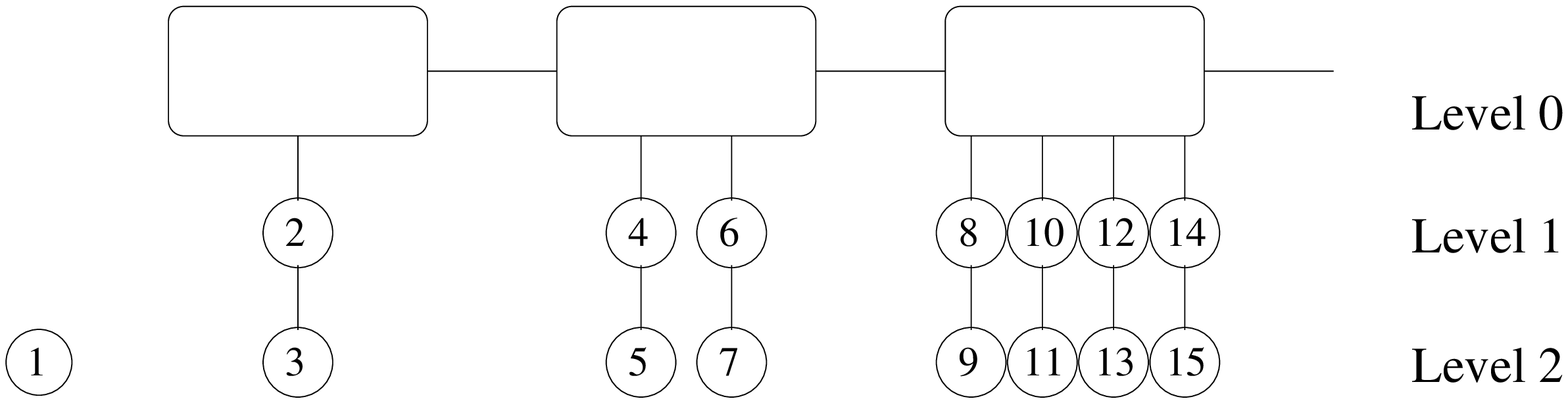}
\end{center}
\caption{The initial portion of $\mathcal{T}_{s}$ with $s=0$.}
\label{fg1}
\end{figure}

Figure \ref{fg1} shows the initial portion of $\mathcal{T}_{s}$. To create $\mathcal{T}_{s}$
first we join an infinite chain of nodes $\{s_i\}_{i=1}^{\infty}$, with node
$s_{i+1}$ connected to node $s_i$ for $i = 1,2,\ldots$. We say that
these nodes, which we distinguish by the term \emph{super-nodes},
are at level 0 in the graph. Each super-node $s_i$ is the root of a
subtree of $\mathcal{T}_{s}$ and has $2^{i-1}$ children, or level 1 nodes. Each level 1 node
has one child, so a level 2 node. Level 2 nodes have no children,
hence these nodes are called \emph{leaves}. Denote by $\mathcal{S}_i$, $i=1,2,\ldots$
the rooted tree consisting of the root $s_i$ and all its descendants
(level 1 and level 2 nodes). Finally we complete $\mathcal{T}_{s}$ by
adding an isolated node $I$. Note that $I$ is itself a leaf of a
trivial tree, but is not of level 0, 1 or 2. As usual, for a graph
$\mathcal{G}$, let $V(\mathcal{G})$ denote the set of nodes
(vertices) of $\mathcal{G}$. Call the subgraph of $\mathcal{T}_{s}$
induced by $\{I\}\cup\bigcup_{k=1}^i V(\mathcal{S}_k)$ the $i$th
\emph{complete partial graph}.

\begin{figure}[htbp]
\begin{center}
\epsfxsize=12cm\epsffile{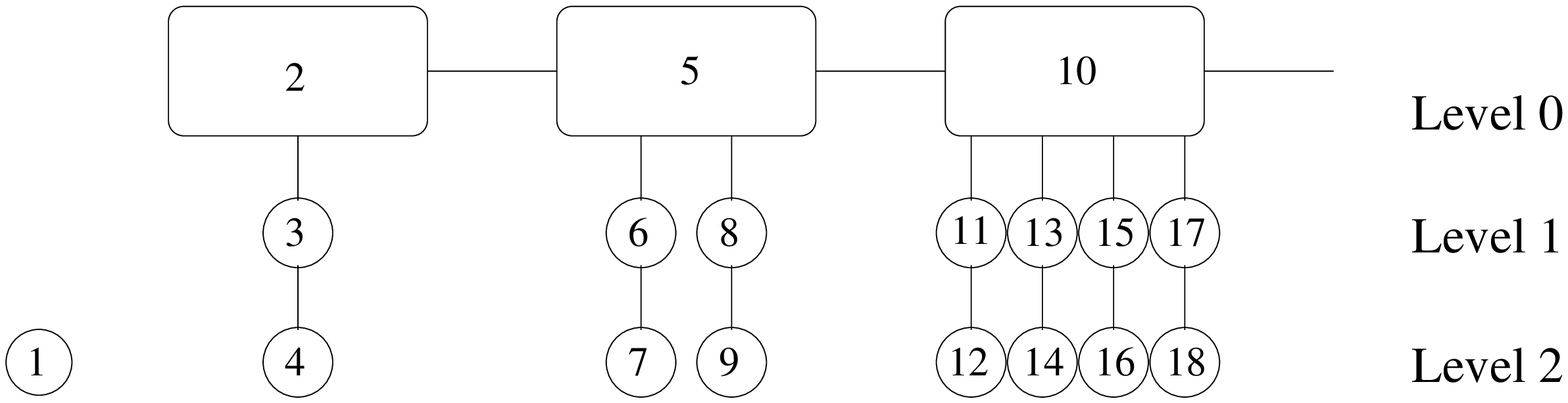}
\end{center}
\caption{The initial portion of $\mathcal{T}_{s}$ with $s=1$.}
\label{fg2}
\end{figure}

\begin{figure}[htbp]
\begin{center}
\epsfxsize=12cm\epsffile{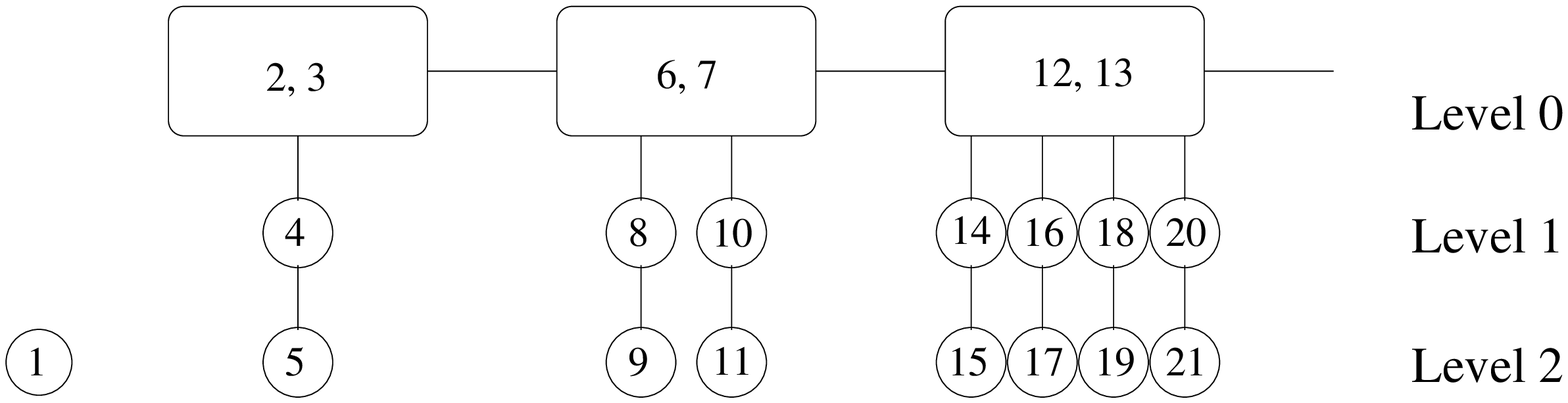}
\end{center}
\caption{The initial portion of $\mathcal{T}_{s}$ with $s=2$.}
\label{fg3}
\end{figure}

We now label each node with the integers $1,2,3,\ldots $ in the
following way. Label $I$, the isolated node, 1. All super-nodes,
starting with $s_1$, have the smallest $s$ consecutive labels not
yet used, where $s$ is the parameter of the  sequence $h_s(n)$ that
we are considering. Note that if $s=0$ the super-nodes have no
labels. So for $s \geq 1$ the first super-node $s_1$ is labelled $2,
3, \ldots, s+1$, which are the next $s$ consecutive integers. Label
each of the level 1 and level 2 nodes of the subtree $\mathcal{S}_1$
with a single integer, in increasing order first from level 1 to
level 2 nodes and then from left to right. Do the same thing for
$\mathcal{S}_2$, $\mathcal{S}_3$ and so on. See Figures \ref{fg1},
\ref{fg2} and \ref{fg3} where we show the initial portion of
$\mathcal{T}_{s}$ for $s=0,1,2$, respectively.

Denote by $v_s(n)$ the node with label $n$. Let $\mathcal{T}_s(n)$
be the subgraph of $\mathcal{T}_{s}$ induced by $\{v_s(i):1\leq i \leq n
\}$. Observe that if $v_s(n)=v_s(n-1)=s_i$ for some $i$, we have
$\mathcal{T}_s(n)=\mathcal{T}_s(n-1)$. Note that the isolated node
$I$ is a leaf of $\mathcal{T}_s(n)$ for all $n=1,2,\ldots$. Let
$a_s(n)$ be the number of leaves of $\mathcal{T}$ in
$\mathcal{T}_s(n)$. That is, $a_s(n)$ equals the number of level 2
nodes in $\mathcal{T}_s(n)$ plus 1 (for the isolated node $I$).
Finally, define $d_s(n)$ to be 1 if $v_s(n)$ is a leaf of
$\mathcal{T}_{s}$ and to be 0 otherwise. Then it follows directly from
these definitions that

\begin{equation}   \label{eqnAD}
a_s(n)=\sum_{k=1}^n d_s(k).
\end{equation}

It is immediate from (\ref{eqnAD}) that the sequence
$\{a_s(n)\}_{n=1}^{\infty}$ is slowly growing. In addition

\begin{lemma}  \label{lmAn}
The  sequence $\{a_s(n)\}_{n=1}^{\infty}$ satisfies the following properties:\\
\emph{(a)} Let $s>0$ and $r \geq 1$. Suppose $n$ is the $r$th label in a super-node. Then for $-r\leq i\leq
s+1-r$, $ a_s(n-r-1)+1=a_s(n+i)=a_s(n-r+s+2)-1.$\\
\emph{(b)} If $v_s(n)$ is a level 1 node, then $ a_s(n-1)=a_s(n)=a_s(n+1)-1.$\\
\emph{(c)} If $v_s(n)$ is a level 2 node, then $ a_s(n-1)+1=a_s(n)=a_s(n+1).$
\end{lemma}

\begin{proof}
(a) Observe that if $n$ is the $r$-th label in a super-node, then $n-r+1$ and $n-r+s$ are the first and last labels respectively in this super-node. Further, $v_s(n-r)$ is the last node in the immediately preceding complete partial graph, so is a leaf. Thus, $v_s(n-r-1)$ is the parent of this leaf and is a child of a super-node. It follows from (\ref{eqnAD}) and our definitions that $a_s(n-r-1)+1$ counts the number of leaves in that preceding complete subtree.

Since $-r \leq i \leq s+1-r$, then $n-r \leq n+i \leq n-r+s+1$. All the values of $n+i$ lie between the label of the last leaf in the immediately preceding complete subtree and the label of the first child of this super-node. Thus,
$$
a_s(n-r-1)+1=a_s(n+i)=a_s(n-r+s+2)-1.
$$

Both (b) and (c) are immediate from the definitions of level 1 and level 2 and (\ref{eqnAD}).
\end{proof}

\begin{lemma} \label{lmCompPtlTr}
For $i \geq 1$, $1\leq r \leq 2^{i-1}$,\\
\emph{(a)} If $v_s(n)$ is the super-node $s_i$ (written compactly $v_s(n)=s_i$), then $a_s(n)=2^{i-1}$.\\
\emph{(b)} If $v_s(n)$ is a level 1 node, the $r$th child of $s_i$, then
$a_s(n)=2^{i-1}+r-1$.\\
\emph{(c)} If $v_s(n)$ is a level 2 node, a leaf whose parent is the $r$th
child of $s_i$, then $a_s(n)=2^{i-1}+r$.

\end{lemma}

\begin{proof}
If $v_s(n)=s_i$, $a_s(n)$ counts the number of leaves in the
$(i-1)$st complete partial graph. So
\begin{equation}
a_s(n)=1+\sum_{k=1}^{i-1} 2^{k-1}=2^{i-1}.
\end{equation}

The other two cases follow readily.
\end{proof}

\begin{figure}[htbp]
\begin{center}
\epsfxsize=10cm\epsffile{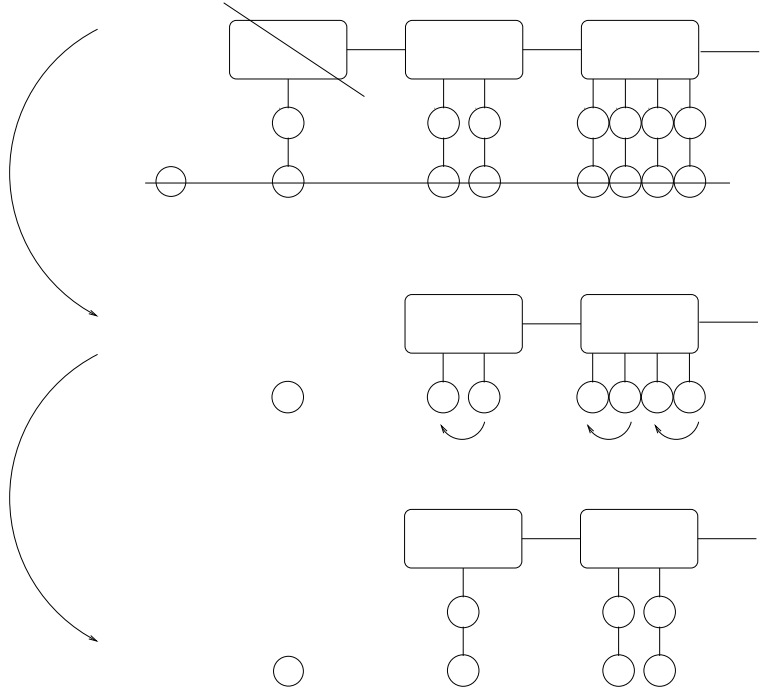}
\caption{The chopping process $P$} \label{fgChop}
\end{center}
\end{figure}

For any $\mathcal{T}_s(n)$, we define a chopping process $P$ that
removes and rearranges nodes in the following way (see Figure
\ref{fgChop}): begin by removing the first super-node $s_1$ and all
leaves of $\mathcal{T}_{s}$ (including the isolated node $I$). Then let
the only child of $s_1$ (which is not removed) be the new isolated
node $I$, and every second level 1 node in $\mathcal{S}_i$ ($i>1$)
be the child of its neighboring sibling on the left. Finally,
relabel the tree so obtained in the same way as we label
$\mathcal{T}_{s}$.

\begin{figure}[htbp]
\begin{center}
\epsfxsize=12cm\epsffile{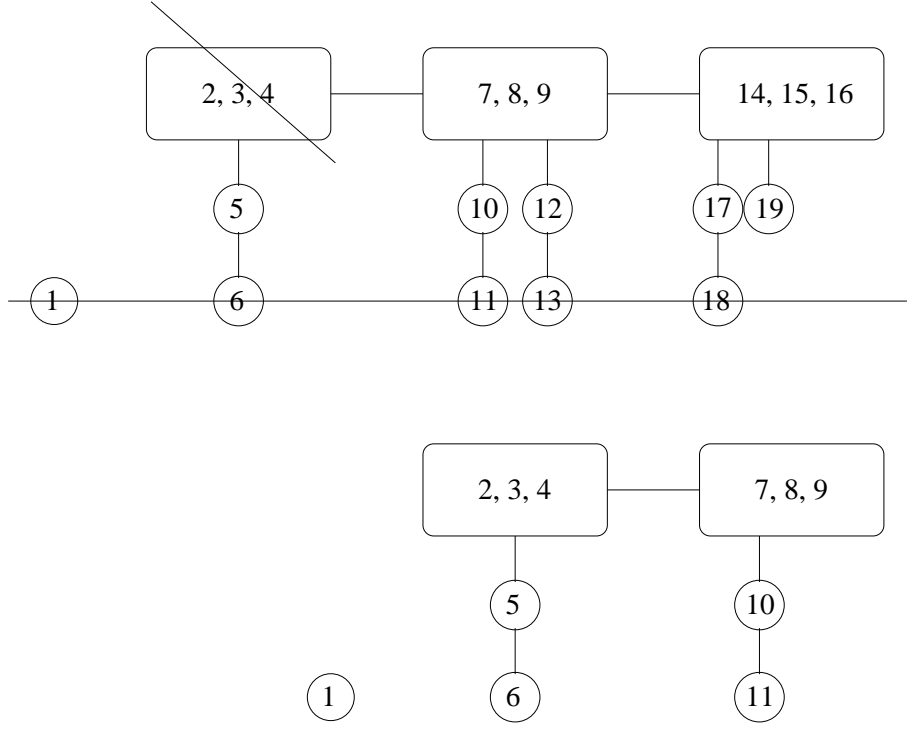}
\caption{The chopping process $P$ reduces $\mathcal{T}_3(19)$ to $\mathcal{T}_3(11)$.}
\label{fgExample}
\end{center}
\end{figure}

We observe several properties of the chopping process
$P$:

\begin{itemize}
\item \textbf{P1.} The tree $\mathcal{T}_{s}$ is invariant under $P$. One still gets $\mathcal{T}_{s}$ after applying $P$ on
$\mathcal{T}_{s}$. See Figure \ref{fgChop}.

\item \textbf{P2.} Apply $P$ on the $i$th complete partial graph, $i>1$. We get the $(i-1)$st complete partial graph.

\item \textbf{P3.} For any $n\geq s+4$, apply $P$ on $\mathcal{T}_s(n)$
(the subgraph $\mathcal{T}_s(n)$ strictly contains the first
complete subgraph). The resulting reduced graph is
$\mathcal{T}_s(n-s-a_s(n))$, which has $s+a_s(n)$ fewer labels (note
that $n-s-a_s(n) > 0$ by construction). For $i>1$,
$\mathcal{T}_s(n)$ contains $\mathcal{S}_i$ if and only if
$\mathcal{T}_s(n-s-a_s(n))$ contains $\mathcal{S}_{i-1}$.  See
Figure \ref{fgExample} where we apply $P$ to $\mathcal{T}_3(19)$ and
get $\mathcal{T}_3(19-3-a_3(19))=\mathcal{T}_3(11)$.

\end{itemize}

Now we state our main theorem:

\begin{theorem}  \label{thmTrHs}
The sequence $h_s(n)$ defined by (\ref{eqnHs}) and the given initial
conditions counts the number of leaves in $\mathcal{T}_s(n)$. More
precisely, for every positive integer $n$, $h_s(n)=a_s(n)$.
\end{theorem}

\begin{proof}
We proceed by induction on $n$ for fixed $s$. We first show that the
assertion holds for the initial values $n=1,2,\ldots , 2s+7$. These
comprise all the labels in the first two rooted trees
$\mathcal{S}_1$ and $\mathcal{S}_2$ (together with $I$).

For $n=1,2,\ldots , s+2$, $h_s(n)=1$ by the initial conditions
assumed for (\ref{eqnHs}), while $a_s(n)=1$ because the isolated
node $I$ has label 1 and counts as a leaf, the super-node $s_1$ has
$s$ labels and its only child is labeled $s+2$. For $n=s+3$
$a_s(n)=2$ because the node labeled $s+3$ is the first (and only)
leaf in $\mathcal{S}_1$. Further, for $n=s+4, \ldots , 2s + 4$,
$a_s(n)=2$ because the super-node $s_2$ contains $s$ labels starting
with $s+4$ and the node labeled $2s+4$ is the first child of $s_2$.
But for the range of values $n=s+3, \ldots , 2s + 4$ we can apply
induction to the recursion for $h_s(n)$ to get
$h_s(n)=h_s(n-s-h_s(n-1))+h_s(n-2-s-h_s(n-3))=1+1=2=a_s(n)$, as
desired. The next two labels $n=2s+5$ and $2s+6$ are assigned to the
first leaf in $\mathcal{S}_2$ and the second child of $s_2$,
respectively, so it follows that $a_s(2s+5)=a_s(2s+6)=3$. But for
these values of $n$, $h_s(n)=3$ by (\ref{eqnHs}). Finally, the label
$2s+7$ is assigned to the second (and last) leaf in $\mathcal{S}_2$,
so once again applying the above recursion we get
$h_s(2s+7)=4=a_s(2s+7)$. See Figure \ref{fg3}.

We prove the assertion for $n$ assuming the assertion is true for
the first $n-1$ entries, where $v_s(n)$ is not a node of the 2nd
complete partial graph (i.e., $n > 2s+7$).
\\
\\
\emph{Case 1}: Node $v_s(n)$ is a level 0 node, i.e., a super-node. Assume
$v_s(n)=s_i, i\geq 3$. By Lemma \ref{lmCompPtlTr}(a), we have
$a_s(n)=2^{i-1}$. We will prove that $
a_s(n-s-a_s(n-1))=a_s(n-2-s-a_s(n-3))=2^{i-2}$, from which it
follows that $a_s(n)=2^{i-1}=a_s(n-s-a_s(n-1))+a_s(n-2-s-a_s(n-3))$.

By the induction assumption $a_s(n-1)=h_s(n-1)$ and
$a_s(n-3)=h_s(n-3)$. Thus, $a_s(n-s-a_s(n-1))=a_s(n-s-h_s(n-1))$ and
$a_s(n-2-s-a_s(n-3))=a_s(n-2-s-h_s(n-3))$. Apply the induction
assumption once again to conclude that
$a_s(n-s-h_s(n-1))=h_s(n-s-h_s(n-1))$ and
$a_s(n-2-s-h_s(n-3))=h_s(n-2-s-h_s(n-3))$. But by (\ref{eqnHs})
$h_s(n-s-h_s(n-1))+h_s(n-2-s-h_s(n-3))=h_s(n)$, so we get
$a_s(n)=h_s(n)$, as desired.

To complete Case 1, we need to show that $a_s(n-s-a_s(n-1))=a_s(n-2-s-a_s(n-3))=2^{i-2}$.
Since $v_s(n)= s_i$, we have $a_s(n-1)=a_s(n)$ by Lemma \ref{lmAn}(a).
After the chopping process $P$ on $\mathcal{T}_s(n)$, we get by P3 a new reduced
graph with the last node $v_s(n-s-a_s(n))=s_{i-1}$. By
Lemma \ref{lmCompPtlTr},
$a_s(n-s-a_s(n))=2^{i-2}$, so $a_s(n-s-a_s(n-1))=2^{i-2}$.

If $v_s(n-2)=s_i$, let $n':=n-2$. Then
$a_s(n-2-s-a_s(n-3))=a_s(n'-s-a_s(n'-1))=2^{i-2}$. Otherwise,
$v_s(n-2)$ is a level 2 or level 1 node, so by Lemma \ref{lmAn},
$a_s(n-3)=a_s(n-1)-1$. Thus
$a_s(n-2-s-a_s(n-3))=a_s(n-2-s-a_s(n-1)+1)=a_s(n-1-s-a_s(n-1))$. But
we have already shown that
$v_s(n-s-a_s(n))=v_s(n-s-a_s(n-1))=s_{i-1}$, so by Lemmas \ref{lmAn}
and \ref{lmCompPtlTr},
$a_s(n-1-s-a_s(n-1))=a_s(n-s-a_s(n-1))=2^{i-2}$. This completes Case
1.
\\
\\
\emph{Case 2}: Node $v_s(n)$ is a level 1 node. Assume $v_s(n)$ is a node
of subtree $\mathcal{S}_i$ and the $r$th child of $s_i$, $i \geq 3$.
By Lemma \ref{lmCompPtlTr}(b), we have $a_s(n)=2^{i-1}+r-1$. For
$x\in\R$, let $\lfloor x \rfloor$ be the floor function of $x$. We
will prove that $a_s(n-s-a_s(n-1))=2^{i-2}+\lfloor \frac r2 \rfloor$
and $a_s(n-2-s-a_s(n-3))=2^{i-2}+\lfloor \frac {r-1}{2} \rfloor$,
from which it follows that
$a_s(n)=2^{i-1}+r-1=a_s(n-s-a_s(n-1))+a_s(n-2-s-a_s(n-3))=h_s(n-s-h_s(n-1))+h_s(n-2-s-h_s(n-3))=h_s(n)$,
from the induction assumption as in Case 1.

Since $v_s(n)$ is a level 1 node, we have $a_s(n-1)=a_s(n)$ by Lemma
\ref{lmAn}(b). After the chopping process $P$ on $\mathcal{T}_s(n)$, we get by P3 a
reduced graph with the last node $v_s(n-s-a_s(n))$, the $(\lfloor
\frac r2 \rfloor+1)$st child ($r$ odd) of $s_{i-1}$ or the leaf whose parent
is the $\lfloor \frac r2 \rfloor$th child ($r$ even) of $s_{i-1}$.
By Lemma \ref{lmCompPtlTr},
$a_s(n-s-a_s(n))=2^{i-2}+\lfloor \frac r2 \rfloor$, so $a_s(n-s-a_s(n-1))=2^{i-2}+\lfloor \frac r2 \rfloor$.

If $v_s(n-2)$ is also a child of $s_i$, then it is the $(r-1)$st
child. Let $n':=n-2$. Then
$a_s(n-2-s-a_s(n-3))=a_s(n'-s-a_s(n'-1))=2^{i-2}+\lfloor
\frac{r-1}{2} \rfloor$. Otherwise, $r=1$. We have two subcases:

(a) When $s>1$, then $a_s(n)=a_s(n-3)$. After the chopping process
$P$, by P3 the last node of the reduced graph is
$v_s(n-s-a_s(n-3))$, the first child of $s_{i-1}$. So
$v_s(n-2-s-a_s(n-3))=s_{i-1}$. Thus by Lemma \ref{lmCompPtlTr},
$a_s(n-2-s-a_s(n-3))=2^{i-2}=2^{i-2}+\lfloor \frac {r-1}{2}
\rfloor.$

(b) When $s=0$ or 1, then $a_s(n)=a_s(n-3)+1$. After the chopping
process $P$, by P3 the last node of the reduced graph is
$v_s(n-s-a_s(n))$, the first child of $s_{i-1}$. Thus by Lemma
\ref{lmAn} and Lemma \ref{lmCompPtlTr},
$a_s(n-2-s-a_s(n-3))=a_s(n-s-a_s(n)-1)=a_s(n-s-a_s(n))=2^{i-2}=2^{i-2}+\lfloor
\frac {r-1}{2} \rfloor.$
\\
\\
\emph{Case 3}: $v_s(n)$ is a level 2 node, i.e., a leaf. Assume the
parent of $v_s(n)$ is the $r$th child of $s_i$, $i \geq 3$. By Lemma
\ref{lmCompPtlTr}(c), we have $a_s(n)=2^{i-1}+r$. We will prove that
$a_s(n-s-a_s(n-1))=2^{i-2}+\lfloor \frac {r+1}{2} \rfloor$ and
$a_s(n-2-s-a_s(n-3))=2^{i-2}+\lfloor \frac {r}{2} \rfloor$, from
which it follows that
$a_s(n)=2^{i-1}+r=a_s(n-s-a_s(n-1))+a_s(n-2-s-a_s(n-3))=h_s(n-s-h_s(n-1))+h_s(n-2-s-h_s(n-3))=h_s(n)$,
 from the induction assumption as in Case 1.

Since $v_s(n)$ is a level 2 node, we have $a_s(n-1)=a_s(n)-1$ by
Lemma \ref{lmAn}(c). After the chopping process $P$ on $\mathcal{T}_s(n)$, we
get by P3 a reduced graph with the last node $a_s(n-s-a_s(n-1))$ the
$(\lfloor \frac r2 \rfloor+1)$st child of $s_{i-1}$ or the leaf
whose parent is the $\lfloor \frac r2 \rfloor$th child of $s_{i-1}$. By Lemma \ref{lmCompPtlTr},
$a_s(n-s-a_s(n)+1)=2^{i-2}+\lfloor \frac {r+1}{2}\rfloor$, so $a_s(n-s-a_s(n-1))=2^{i-2}+\lfloor \frac {r+1}{2}\rfloor$.

If $v_s(n-2)$ is also a node of $\mathcal{S}_i$, then it is the
child of the $(r-1)$st child of $s_i$. Let $n':=n-2$. Then
$a_s(n-2-s-a_s(n-3))=a_s(n'-s-a_s(n'-1))=2^{i-2}+\lfloor \frac
{r}{2} \rfloor$. Otherwise, $r=1$. We have two subcases:

(a) When $s>0$, then $a_s(n)=a_s(n-3)$. After the chopping process
$P$, by P3 the last node of the reduced graph is $v_s(n-s-a_s(n))$,
the first child of $s_{i-1}$. So $v_s(n-s-a_s(n)-1)=s_{i-1}$. Thus
by Lemma \ref{lmAn} and Lemma \ref{lmCompPtlTr},
$a_s(n-s-a_s(n-3)-2)=a_s(n-s-a_s(n)-2)=a_s(n-s-a_s(n)-1)=2^{i-2}=2^{i-2}+\lfloor
\frac {r}{2} \rfloor.$

(b) When $s=0$, then $a_s(n)=a_s(n-3)+1$. After the chopping process
$P$, by P3 the last node of the reduced tree is $v_s(n-s-a_s(n))$,
the first child of $s_{i-1}$. Thus by Lemma \ref{lmAn} and Lemma
\ref{lmCompPtlTr},
$a_s(n-2-s-a_s(n-3))=a_s(n-s-a_s(n)-1)=a_s(n-s-a_s(n))=2^{i-2}=2^{i-2}+\lfloor
\frac {r}{2} \rfloor.$

This completes the proof.

\end{proof}

 From Theorem \ref{thmTrHs}, Lemma \ref{lmAn} and the definition of
$a_s(n)$ the following result is immediate .

\begin{corollary} \label{corHsIncs}
For every $s \geq 0$, the sequence $h_s(n)$ is slowly growing. The
sequence hits every positive integer twice except for the powers of
2; the number of repetitions of $2^r, r\geq 0$, in the sequence
$h_s(n)$ is $s+2$.
\end{corollary}

For $s=0$ it is evident from the definition of $a_0(n)$ and Theorem
\ref{thmTrHs} that $h_0(n)=\lceil \frac {n}{2} \rceil$, which
provides a simple closed form for $h_0(n)$. This formula is
immediate from a casual inspection of Table \ref{tab1}, and also
follows directly from Corollary \ref{corHsIncs}.

Observe that Corollary \ref{corHsIncs} implies that $2^{k}+sk$ is
the index of the last term in $h_s(n)$ that equals $2^{k-1}$. The
argument for this is as follows: there are $k$ powers of 2 less than
or equal to $2^{k-1}$, each of which occurs $s+2$ times in $h_s(n)$
while each of the other $2^{k-1}-k$ positive integers less than or
equal to $2^{k-1}$ occurs twice. But since
$2(2^{k-1}-k)+(s+2)k=2^{k}+sk$ we can write $h_s(n)$ in terms of
$h_0(n)$. Since we have a closed form for $h_0(n)$ this provides a
``piecewise'' closed formula for $h_s(n)$ where $s\geq 1$:

\begin{equation}
 h_{s}(n)= \left \{ \begin{array}{ll}
 h_0 (n-sk) & 2^{k}+sk \leq n < 2^{k+1}+sk \\
 2^k & 2^{k+1}+sk \leq n < 2^{k+1}+s(k+1)
  \end{array} \right.
\end{equation}

In a surprising twist it turns out that the sequence $h_1(n)$ also
has an explicit closed form. This sequence, with the additional
initial term $h_1(0) =1$, appears in \cite{sna} as entry A109964
(Bottomley's sequence), though the latter is generated by a very
different recursion that is not self-referencing. To prove that the
two sequences are the same it is enough to show that for $n>0$
Bottomley's sequence satisfies Corollary \ref{corHsIncs}. This
follows from a straightforward induction; we omit the details. Thus,
for $n>0$ $h_1(n)$ satisfies the recursion in \cite{sna} for
A109964(n), namely, $h_1(n)=\lfloor \sqrt{\sum_{j=0}^{n-1} h_1(j)}
\rfloor$.

Further, from \cite{sna} we have that for $n>0$ $A109964(n)=
\lfloor\sqrt{b(n-1)}\rfloor$ where $b(n)$ is defined by the
recursion $b(n)=b(n-1) + \lfloor \sqrt {b(n-1)} \rfloor$, with
$b(0)=1$. But in \cite{gkp} chapter 3, exercise 3.28, the following
closed form for $b(n)$ is attributed to Carl Witty:
$$
b(n-1)=2^{k}+\lfloor (\frac {n-k}{2})^2 \rfloor, \qquad 2^k+k\leq
n<2^{k+1}+k+1.
$$

Thus we have the following closed form for $h_1(n)$:

$$
h_1(n)=\lfloor
\sqrt{2^{k}+\lfloor(\frac{n-k}{2})^{2}\rfloor}\rfloor, \qquad
2^{k}+sk \leq n \leq 2^{k+1}+sk.
$$

For $s>1$ there does not seem to be any comparable closed form for
$h_s(n)$.

Note that it is possible to generate the family of sequences
$h_s(n)$ from a slightly different infinite tree graph and a
modified numbering scheme for the nodes.

\begin{figure}[htbp]
\begin{center}
\epsfxsize=12cm\epsffile{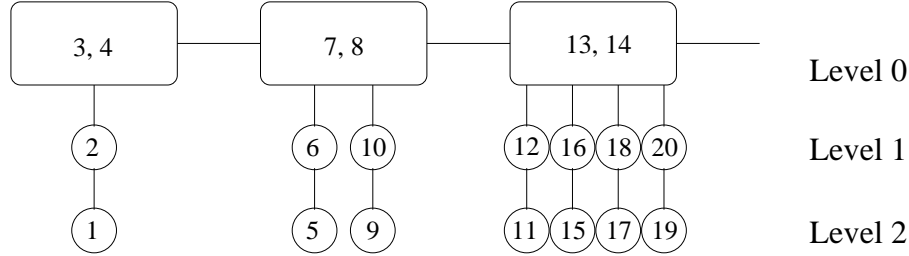}
\end{center}
\caption{The initial portion of $\mathcal{T'}_{s}$ with $s=2$.}
\label{fg6}
\end{figure}

The new tree structure $\mathcal{T'}_s$ is our original
$\mathcal{T}_s$ minus the isolated node $I$. We label each node with
integers 1, 2, 3, $\cdots$ as follows: All subtrees, starting with
$\mathcal{S}_1$, are labelled consecutively with the smallest labels
not yet used. For each subtree $\mathcal{S}_i$, first label the leaf
of the first child of super-node $s_i$ and then label the first
child of $s_i$, each node receiving one label. Next, label the
super-node $s_i$ with $s$ consecutive labels not yet used, where $s$
is the parameter of the sequence $h_s(s)$. Now label each of the
remaining level 2 and level 1 nodes of the subtree $\mathcal{S}_i$
with a single label in increasing order first from level 2 to level
1 nodes and then from left to right. See Figure \ref{fg6}, where we
show the initial portion of $\mathcal{T'}_s$ for $s=2$. If we define
the chopping process in the same way as before, with the obvious
modification since we no longer have an isolated node, then it is
easy to see that the sequence $h_s(n)$ can be interpreted in terms
of the leaves in the this new graph $\mathcal{T'}_s$.

\section{Generating Function} \label{secGenFn}

We compute a generating function for $\{a_s(n)\}_{n=1}^{\infty}$,
which is also the generating function for $\{h_s(n)\}_{n=1}^{\infty}$ by Theorem \ref{thmTrHs}. Define
$\mathcal{A}_s(z):=\sum_{k=1}^{\infty}a_s(k)z^k$ and
$\mathcal{D}_s(z):=\sum_{k=1}^{\infty}d_s(k)z^k$.
By (\ref{eqnAD}), we know that $\mathcal{A}_s(z)=\frac{\mathcal{D}_s(z)}{1-z}$. So it's enough to find $\mathcal{D}_s(z)$.
 From the definition of $d_s(n)$ it is clear that the subsequence of the sequence $d_s(n)$ corresponding to
$\mathcal{S}_i$, $i>0$, is of the form $0^s 01\cdots01$,
with $2^{i-1}$ consecutive $01$'s, preceded by $s$ consecutive 0's, denoted by $0^s$.
Thus $\{d_s(n)\}_{n=1}^{\infty}$ is
$1\vdots 0^s01 \vdots 0^s0101 \vdots 0^s01010101 \vdots \cdots$

For $n>1$, define the concatenation $C_n=C_{n-1}C_{n-1}$, with $C_0=1$ and $C_1=01$.
Let $C_n(k)$ be the $k$th entry of $C_n$. It is obvious that $\{d_s(n)\}_{n=1}^{\infty}$ can be written as
$C_0  0^s C_1 0^s C_2  0^s C_3 0^s\cdots$. Further we have

\begin{lemma}
The length of $C_n$ is $|C_n|=2^n$ and for $n>0$,
$$
\mathcal{C}_n(z):=\sum_{k=1}^{2^n}C_n(k)z^k=\sum_{i=1}^{2^{n-1}}
z^{2i}=\frac{z^2}{1-z^2}(1-z^{2^n}).
$$
\end{lemma}

\begin{proof}
Evidently $|C_0 |=1$ and for $n>0$, $C_n$ has $2^{n-1}$ pairs of 01, so $|C_n|=2^n$. Each 1 occurs as the second entry in every pair, so the result follows.
\end{proof}

\begin{theorem}
The generating function $\mathcal{D}_s(z)$ is equal to
$$
 z+ \frac{z^2}{1-z^2} \sum_{i=1}^{\infty} z^{si+2^i-1}(1-z^{2^i})
$$
\end{theorem}
\begin{proof}
\begin{align*}
\mathcal{D}_s(z)&= z+ \sum_{i=1}^{\infty} z^{si+\sum_{j=0}^{i-1}|C_j|}\mathcal{C}_i(z) \\
              &= z+ \sum_{i=1}^{\infty} z^{si+2^i -1}\frac{z^2}{1-z^2}(1-z^{2^i}) \\
              &= z+ \frac{z^2}{1-z^2} \sum_{i=1}^{\infty} z^{si+2^i-1}(1-z^{2^i})
\end{align*}

\end{proof}

\begin{remark}
We can also rewrite the expression for $\mathcal{D}_s(z)$
\begin{align*}
\mathcal{D}_s(z)&= z+ \frac{z^2}{1-z^2} \sum_{i=1}^{\infty} z^{si+2^i-1}(1-z^{2^i})\\
                                    & =  z+ \frac{z^2}{1-z^2}z^{s+2^0}[1-z^{2^1}+z^{s+2^1}[1-z^{2^2}+z^{s+2^2}[1-z^{2^3}+z^{s+2^3}[\cdots\\
                                    &=\frac{z^2}{1-z^2}[\frac1z-z^{2^0}+z^{s+2^0}[1-z^{2^1}+z^{s+2^1}[1-z^{2^2}+z^{s+2^2}[1-z^{2^3}+z^{s+2^3}\cdots
\end{align*}
\end{remark}

\begin{corollary}
The generating function $\mathcal{A}_s(z)$ is equal to
$$
 \frac{z}{1-z}+ \frac{z^2}{(1-z^2)(1-z)} \sum_{i=1}^{\infty} z^{si+2^i-1}(1-z^{2^i})
$$
\end{corollary}

\begin{proof}
By the proposition and $\mathcal{A}_s(z)=\frac{\mathcal{D}_s(z)}{1-z}$.
\end{proof}

\section{Concluding Remark}  \label{secRemark}
It is equally natural to investigate the meta-Fibonacci recursion that derives from the odd values of the index on the right hand side of (\ref{equTsk}), that is,
\begin{equation}  \label{eqnGsk}
 g_{s,k}(n)=\sum_{j=0}^{k-1} g_{s,k}(n- (2j+1) -s - g_{s,k}(n-(2j+2))), \qquad n>s+2k.
\end{equation}
For general $k$ and $s$, the recursion (\ref{eqnGsk}) produces a
much richer collection of infinite sequences (that is, sequences
that do not die) for various sets of initial conditions. In many
instances these sequences exhibit considerably more complicated
behavior than $h_s(n)$. To date we have not been able to discover a
combinatorial interpretation for any of these sequences, whose
behavior will be the subject of a future communication.

\section*{Acknowledgements}\label{ack}
The authors are grateful to Rafel Drabek, Noah Cai and Junmin Park
for research assistance at various early stages of this work, and to
an anonymous referee for several helpful observations.

\bigskip
\hrule
\bigskip

\noindent 2000 {\it Mathematics Subject Classification:} Primary 05A15. Secondary 11B37; 11B39.

\noindent \emph{Keywords:}   meta-Fibonacci recursion; Hofstadter sequence.

\bigskip
\hrule
\bigskip

\noindent (Concerned with sequences \seqnum{A008619}, \seqnum{A109964}.)

\bigskip
\hrule
\bigskip

\noindent Return to \htmladdnormallink{Journal of Integer Sequences
home page}{http://www.cs.uwaterloo.ca/journals/JIS/}. \vskip .1in

\end{document}